\documentclass[11pt]{amsart}

\usepackage{pstricks,graphicx,pst-node,rotating,amsmath,amssymb}
\newtheorem{theorem}{Theorem}[section]

\newtheorem{proposition}[theorem]{Proposition}
\newtheorem{corollary}[theorem]{Corollary}

\newtheorem{definition}[theorem]{Definition}

\newtheorem{remark}[theorem]{Remark}

\newcommand{\reffig}[1]{Figure \ref{fig:#1}}
\newcommand{\refeq}[1]{equation (\ref{eqn:#1})}

\newcommand{\refthm}[1]{Theorem \ref{thm:#1}}

\newcommand{\refprop}[1]{Proposition \ref{prop:#1}}

\newcommand{\refdef}[1]{Definition \ref{defn:#1}}

\newcommand{\refsec}[1]{Section \ref{sec:#1}}

\newcommand{\Sn}{\mathcal{S}_n}

\newcommand{\vp}{\varphi}
\newcommand{\SN}{S \setminus n}

\newcommand{\wJ}{\widetilde{J}}

\newcommand{\row}{\mathrm{row}}
\newcommand{\col}{\mathrm{col}}

\newcommand{\rowd}{\succ}
\newlength{\hsp}
\newlength{\vsp}
\newlength{\vspi}
\newlength\cellsize 
\savebox2{%
\begin{picture}(12,12)
\put(0,0){\line(1,0){12}}
\put(0,0){\line(0,1){12}}
\put(12,0){\line(0,1){12}}
\put(0,12){\line(1,0){12}}
\end{picture}}
\newcommand\cellify[1]{\def\thearg{#1}\def\nothing{}%
\ifx\thearg\nothing
\vrule width0pt height\cellsize depth0pt\else
\hbox to 0pt{\usebox2\hss}\fi%
\vbox to 12\unitlength{
\vss
\hbox to 12\unitlength{\hss$#1$\hss}
\vss}}
\newcommand\tableau[1]{\vtop{\let\\=\cr
\setlength\baselineskip{-12000pt}
\setlength\lineskiplimit{12000pt}
\setlength\lineskip{0pt}
\halign{&\cellify{##}\cr#1\crcr}}}

\newlength\smcellsize 
\savebox3{%
\begin{picture}(8,8)
\put(0,0){\line(1,0){8}}
\put(0,0){\line(0,1){8}}
\put(8,0){\line(0,1){8}}
\put(0,8){\line(1,0){8}}
\end{picture}}
\newcommand\smcellify[1]{\def\thearg{#1}\def\nothing{}%
\ifx\thearg\nothing
\vrule width0pt height\smcellsize depth0pt\else
\hbox to 0pt{\usebox3\hss}\fi%
\vbox to 8\unitlength{
\vss
\hbox to 8\unitlength{\hss$#1$\hss}
\vss}}
\newcommand\smtableau[1]{\vtop{\let\\=\cr
\setlength\baselineskip{-10000pt}
\setlength\lineskiplimit{10000pt}
\setlength\lineskip{0pt}
\halign{&\smcellify{##}\cr#1\crcr}}}

\newcommand{\tabloidXxx}[3]{%
  \psset{xunit=1\cellsize}
  \psset{yunit=1\cellsize}
  \pspicture(0,0)(2,2)
  \pspolygon(0,0)(0,2)(1,2)(1,1)(2,1)(2,0)
  \psline(0,1)(2,1)
  \rput(0.5,1.5){$#1$}
  \rput(0.5,0.5){$#2$} \rput(1.5,0.5){$#3$}
  \endpspicture}
\newcommand{\tabXxx}[3]{%
  \psset{xunit=1\cellsize}
  \psset{yunit=1\cellsize}
  \pspicture(0,0)(2,2)
  \pspolygon(0,0)(0,2)(1,2)(1,1)(2,1)(2,0)
  \psline(0,1)(2,1)  \psline(1,0)(1,1)
  \rput(0.5,1.5){$#1$}
  \rput(0.5,0.5){$#2$} \rput(1.5,0.5){$#3$}
  \endpspicture}

\begin{document}

%%%%%%%%%%%%%%%%%%%%%%%%%%%%%%%%%%%%%%%%%%%%%%%%%%%%%%%%%%%%
%  TITLE PAGE information
%%%%%%%%%%%%%%%%%%%%%%%%%%%%%%%%%%%%%%%%%%%%%%%%%%%%%%%%%%%%

%     [Short Title]{Full Title}
\title[Kicking the $n!$ Theorem]{A kicking basis for the two column
  Garsia-Haiman modules}

\author[S. Assaf]{Sami Assaf} %
\address{Department of Mathematics, Massachusetts Institute of
  Technology, 77 Massachusetts Avenue, Cambridge, MA 02139-4307}
%\curraddr{}
\email{sassaf@math.mit.edu}
%\thanks{}

\author[A. Garsia]{Adriano Garsia} %
\address{Department of Mathematics, University of California, San
  Diego, 9500 Gilman Drive, La Jolla, CA 92093-0112}
%\curraddr{}
\email{garsia@ucsd.edu}
%\thanks{}

%    General info
 % Actions of groups on commutative rings; invariant theory


\begin{abstract}
  In the early 1990s, Garsia and Haiman conjectured that the dimension
  of the Garsia-Haiman module $R_{\mu}$ is $n!$, and they showed that
  the resolution of this conjecture implies the Macdonald Positivity
  Conjecture. Haiman proved these conjectures in 2001 using algebraic
  geometry, but the question remains to find an explicit basis for
  $R_{\mu}$ which would give a simple proof of the dimension. Using
  the theory of Orbit Harmonics developed by Garsia and Haiman, we
  present a "kicking basis" for $R_{\mu}$ when $\mu$ has two
  columns. 
\end{abstract}

\maketitle

%%%%%%%%%%%%%%%%%%%%%%%%%%%%%%%%%%%%%%%%%%%%%%%%%%%%%%%%%%%%
\section{Introduction}
%%%%%%%%%%%%%%%%%%%%%%%%%%%%%%%%%%%%%%%%%%%%%%%%%%%%%%%%%%%%
\label{sec:intro}

In 1988, Macdonald \cite{Macdonald1988} found a remarkable new basis
of symmetric functions in two parameters which specializes to Schur
functions, complete homogeneous, elementary and monomial symmetric
functions and Hall-Littlewood functions, among others. With an
appropriate analog of the Hall inner product, the transformed
Macdonald polynomials $\widetilde{H}_{\mu}(Z;q,t)$ are uniquely
characterized by certain triangularity and orthogonality conditions,
from which their symmetry follows. The {\em Kostka-Macdonald}
polynomials, $\widetilde{K}_{\lambda\mu}(q,t)$, are defined by
$$
\widetilde{H}_{\mu}(Z;q,t) = \sum_{\lambda}
\widetilde{K}_{\lambda\mu}(q,t) s_{\lambda}(Z) .
$$
The Macdonald Positivity Conjecture states that
$\widetilde{K}_{\lambda\mu} (q,t) \in \mathbb{N}[q,t]$.

In 1993, Garsia and Haiman \cite{GaHa1993} conjectured that the
transformed Macdonald polynomials could be realized as the bigraded
characters for a diagonal action of $S_n$ on two sets of
variables. Moreover, they were able to show that knowing the dimension
of this module is enough to determine its character. Therefore the $n!$
Conjecture, which states that the dimension of the Garsia-Haiman
module is $n!$, implies the Macdonald Positivity Conjecture.

By analyzing the algebraic geometry of the Hilbert scheme of $n$
points in the plane, Haiman \cite{Haiman2001} was able to prove the
$n!$ Conjecture and consequently establish Macdonald
Positivity. However, it remains an important open problem in the
theory of Macdonald polynomials to prove the $n!$ Theorem directly by
finding an explicit basis for the module. After reviewing Macdonald
polynomials and the Garsia-Haiman modules in \refsec{mods}, we give an
explicit basis for the Garsia-Haiman modules indexed by a partition
with at most two columns in \refsec{2col}. A new basis for hooks is
also given in \refsec{hooks}.

%%%%%%%%%%%%%%%%%%%%%%%%%%%%%%%%%%%%%%%%%%%%%%%%%%%%%%%%%%%%
\section{Macdonald polynomials and graded $\Sn$-modules}
%%%%%%%%%%%%%%%%%%%%%%%%%%%%%%%%%%%%%%%%%%%%%%%%%%%%%%%%%%%%
\label{sec:mods}

We assume the definitions and notations from \cite{Macdonald1995} of
partitions and the classical bases for symmetric functions. So as to
avoid confusions when defining various modules, we use the alphabet $Z
= z_1,\ldots,z_n$ for symmetric functions. For example, the
\textit{Schur functions} shall be denoted $s_{\lambda}(Z)$. 

%%%%%%%%%%%%%%%%%%%%%%%%%%%%%%%%%%%%%%%%%%%%%%%%%%%%%%%%%%%%
\subsection{Macdonald positivity}
%%%%%%%%%%%%%%%%%%%%%%%%%%%%%%%%%%%%%%%%%%%%%%%%%%%%%%%%%%%%
\label{sec:pos}

Departing slightly from Macdonald's convention of defining
$P_{\mu}(Z;q,t)$ \cite{Macdonald1988}, we instead use the transformed
Macdonald polynomials $\widetilde{H}_{\mu}(Z;q,t)$ as presented in
\cite{GaHa1993}.

\begin{definition}
  The transformed Macdonald polynomials $\widetilde{H}_{\mu}(Z;q,t)$
  are the unique functions satisfying the following triangularity and
  orthogonality conditions:
  \begin{itemize}
    \item[(i)] $\widetilde{H}_{\mu}(Z;q,t) \in \mathbb{Q}(q,t) \{
    s_{\lambda}[Z/(1-q)] \ : \ \lambda \geq \mu\}$;
    \item[(ii)] $\widetilde{H}_{\mu}(Z;q,t) \in \mathbb{Q}(q,t) \{
    s_{\lambda}[Z/(1-t)] \ : \ \lambda \geq \mu'\}$;
    \item[(iii)] $\widetilde{H}_{\mu}[1;q,t] = 1$.
  \end{itemize}
  \label{defn:macdefn}
\end{definition}

The square brackets in \refdef{macdefn} stand for {\em plethystic
  substitution}. In short, $s_{\lambda}[A]$ means $s_{\lambda}$
applied as a $\Lambda$-ring operator to the expression $A$, where
$\Lambda$ is the ring of symmetric functions. For a thorough
account of plethysm, see \cite{Haiman1999}. 

The existence of such a family of functions is a theorem, following in
large part from Macdonald's original proof of existence. Once
established, the symmetry of $\widetilde{H}_{\mu}(Z;q,t)$ follows by
definition. Of particular importance are the change of basis
coefficients from the transformed Macdonald polynomials to the Schur
functions, defined by
\begin{equation}
  \widetilde{H}_{\mu}(Z;q,t) = \sum_{\lambda}
  \widetilde{K}_{\lambda,\mu}(q,t) s_{\lambda}(Z) .
\end{equation}
A priori, the $\widetilde{K}_{\lambda,\mu}(q,t)$ are rational
functions in $q$ and $t$ with rational coefficients. 

\begin{theorem}[\cite{Haiman2001}]
  We have $\widetilde{K}_{\lambda,\mu}(q,t) \in \mathbb{N}[q,t]$.
\label{thm:mp}
\end{theorem}

Macdonald originally conjectured \refthm{mp} when he introduced the
polynomials in 1988. The original proof, due to Haiman in 2001,
realizes $\widetilde{H}_{\mu}(Z;q,t)$ as the bigraded character of
certain modules for the diagonal action of $\Sn$ on $\mathbb{Q}[X,Y]$;
see sections \ref{sec:mods-GH} and \ref{sec:mods-OH}. From this it
follows that the character can be written as a sum of irreducible
representations of $\Sn$ with coefficients in $\mathbb{N}[q,t]$. Under
the Frobenius image, these coefficients exactly give
$\widetilde{K}_{\lambda,\mu}(q,t)$. The aim of this paper is to follow
this method of proof until it departs the realm of representation
theory for algebraic geometry.

It is worth noting that there are now two additional proofs of
Macdonald positivity, both of which utilize an expansion of Macdonald
polynomials in terms of LLT polynomials conjectured by Haglund
\cite{Haglund2004} and proved along with Haiman and Loehr
\cite{HHL2005}. Grojnowski and Haiman \cite{GrHa2007} have a proof
using Kazhdan-Lusztig theory and the first author \cite{Assaf2007-2}
has a purely combinatorial proof.

%%%%%%%%%%%%%%%%%%%%%%%%%%%%%%%%%%%%%%%%%%%%%%%%%%%%%%%%%%%%
\subsection{Garsia-Haiman modules}
%%%%%%%%%%%%%%%%%%%%%%%%%%%%%%%%%%%%%%%%%%%%%%%%%%%%%%%%%%%%
\label{sec:mods-GH}

\begin{figure}[ht]
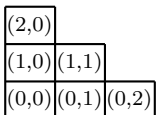

  \begin{center}
    \psset{xunit=4ex}
    \psset{yunit=3ex}
    \pspicture(0,1)(2,3.5)
    \psline(0,0)(3,0) \psline(0,1)(3,1) \psline(0,2)(2,2) \psline(0,3)(1,3)
    \psline(0,0)(0,3) \psline(1,0)(1,3) \psline(2,0)(2,2) \psline(3,0)(3,1)
    \rput(0.5,0.5){$_{(0,0)}$}
    \rput(1.5,0.5){$_{(0,1)}$}
    \rput(0.5,1.5){$_{(1,0)}$}
    \rput(1.5,1.5){$_{(1,1)}$}
    \rput(0.5,2.5){$_{(2,0)}$}
    \rput(2.5,0.5){$_{(0,2)}$}
    \endpspicture
    \caption{\label{fig:coordmu} The coordinates for each cell of $\mu
      = (3,2,1)$.}
  \end{center}    
\end{figure}

To define the modules mentioned in \refsec{pos}, we consider the
diagonal action of the symmetric group $\Sn$ on the polynomial ring
$\mathbb{Q}[X,Y] = \mathbb{Q}[x_1, \ldots, x_n; y_1, \ldots, y_n]$
permuting the $x_i$'s and $y_j$'s simultaneously and identically. Let
the coordinates of the diagram of a partition $\mu$ of $n$ be
$\{(p_1,q_1), \ldots, (p_n,q_n)\}$, where $p$ gives the row coordinate
and $q$ the column coordinate indexed from zero; see \reffig{coordmu}.

Define the polynomial $\Delta_{\mu} \in \mathbb{Q}[X,Y]$ by
\begin{equation}
  \Delta_{\mu}(X,Y) = \det \left(
      \begin{array}{cccc}
        x_{1}^{p_1} y_{1}^{q_1} & x_{2}^{p_1} y_{2}^{q_1} & \cdots &
        x_{n}^{p_1} y_{n}^{q_1} \\[.5\vsp]
        x_{1}^{p_2} y_{1}^{q_2} & x_{2}^{p_2} y_{2}^{q_2} & \cdots &
        x_{n}^{p_2} y_{n}^{q_2} \\
        \vdots & \vdots & & \vdots \\
        x_{1}^{p_n} y_{1}^{q_n} & x_{2}^{p_n} y_{2}^{q_n} & \cdots &
        x_{n}^{p_n} y_{n}^{q_n}
      \end{array} \right).
\end{equation}
Since the bi-exponents are all distinct, $\Delta_{\mu}$ is a non-zero
homogeneous $\Sn$-alternating polynomial with top degree $n(\mu) =
\sum_i (i-1) \mu_i$ in $X$ and $n(\mu')$ in $Y$. Taking $\mu = (1^n)$
or $\mu = (n)$ gives the Vandermonde determinant in $X$ or $Y$,
respectively.

Let $\mathcal{I}_{\mu} \subset \mathbb{Q}[X,Y]$ be the ideal of
polynomials $\vp$ such that
\begin{displaymath}
  \vp(\partial/\partial x_1, \ldots, \partial/\partial
  x_n; \partial/\partial y_1, \ldots, \partial/\partial y_n)
  \Delta_{\mu} = 0.
\end{displaymath}
Clearly this defines an $\Sn$ invariant doubly homogeneous
ideal. Define the \textit{Garsia-Haiman module $\mathcal{H}_{\mu}$} to
be the quotient ring $\mathbb{Q}[X,Y]/\mathcal{I}_{\mu}$ with its
natural structure of a doubly graded $\Sn$-module.

Garsia and Haiman \cite{GaHa1996} proved that if this module has the
correct dimension (the $n!$ Conjecture), then the bi-graded character
is given by the transformed Macdonald polynomial.

\begin{theorem}[\cite{GaHa1996}]
  If $\mathcal{H}_{\mu}$ affords the regular representation of $\Sn$,
  then the bi-graded Frobenius character, given by
  \begin{displaymath}
    \mathrm{Frob}_{\mathcal{H}_{\mu}} (Z; q,t) = \sum_{i,j} t^i q^j
    \psi\left( (\mathcal{H}_{\mu})_{i,j} \right),
  \end{displaymath}
  where $\psi$ is the usual Frobenius map sending the Specht module
  $S^{\lambda}$ to the Schur function $s_{\lambda}$, is equal to the
  transformed Macdonald polynomials $\widetilde{H}_{\mu} (Z; q,t)$. In
  particular, $\widetilde{K}_{\lambda,\mu}(q,t) \in \mathbb{N}[q,t]$.
\label{thm:n2m}
\end{theorem}

The following theorem is the famed $n!$ Conjecture of Garsia and
Haiman \cite{GaHa1993}, proved by Haiman \cite{Haiman2001} in 2001.

\begin{theorem}[\cite{Haiman2001}] 
  The dimension of $\mathcal{H}_{\mu}$ is $n!$.
\label{thm:n!}
\end{theorem}

By \refthm{n2m}, Haiman's proof of the $n!$ Conjecture provided the
first proof of the Macdonald positivity conjecture. Haiman's proof
analyzes the isospectral Hilbert scheme of $n$ points in a plane,
ultimately showing that it is Cohen-Macaulay (and Gorenstein). As this
proof uses difficult machinery in algebraic geometry, it remains an
important open problem to prove \refthm{n!} directly by finding an
explicit basis for the module $\mathcal{H}_{\mu}$.

%%%%%%%%%%%%%%%%%%%%%%%%%%%%%%%%%%%%%%%%%%%%%%%%%%%%%%%%%%%%
\subsection{Orbit Harmonics}
%%%%%%%%%%%%%%%%%%%%%%%%%%%%%%%%%%%%%%%%%%%%%%%%%%%%%%%%%%%%
\label{sec:mods-OH}

Let $\alpha_1, \ldots, \alpha_n$ and $\beta_1, \ldots, \beta_n$ be
sequences of distinct rational numbers. Let $(p_1,q_1), \ldots,
(p_n,q_n)$ be the coordinates of the cells of $\mu$ taken in some
order, recorded by the standard filling $S$ of $\mu$ given by placing
the entry $i$ in the cell $(p_i,q_i)$. To each $S$, associate the {\em
  orbit point of $S$}, denoted $p_S$, defined by
\begin{equation}
  p_S = (\alpha_{p_1 +1}, \ldots, \alpha_{p_n +1}; \beta_{q_1 +1}, \ldots,
  \beta_{q_n +1}) .
\end{equation}
Here the shift in indices is a notational convenience. For example,
\begin{displaymath}
  p \ \raisebox{-.5\cellsize}{\smtableau{ _5 & _3 \\ _6 & _1 \\ _2 & _4}}
  \ = \ (\alpha_2, \alpha_1, \alpha_3, \alpha_1, \alpha_3, \alpha_2;
  \beta_2, \beta_1, \beta_2, \beta_2, \beta_1, \beta_1) .
\end{displaymath}

Let $\Sn$ act on $\mathbb{Q}^{2n}$ by permuting the first $n$ and
second $n$ coordinates simultaneously and identically. Let $[p_S]$
denote the \textit{regular} orbit of $p_S$ under this action.
Regarding $\mathbb{Q}[X,Y]$ as the coordinate ring of
$\mathbb{Q}^{2n}$, define $\mathcal{J}_{\mu} \subset \mathbb{Q}[X,Y]$
to be the ideal of polynomials vanishing on $[p_S]$. Define the module
$R_{\mu}$ to be the coordinate ring of $[p_S]$, i.e.
$\mathbb{Q}[X,Y]/\mathcal{J}_{\mu}$, with its natural $\Sn$ action.

Since $R_{\mu}$ clearly affords the regular representation, the aim is
to relate this module to $\mathcal{H}_{\mu}$. To do this, construct
the associated graded module $\mathrm{gr}R_{\mu} =
\mathbb{Q}[X,Y]/\mathrm{gr}\mathcal{J}_{\mu}$. Garsia and Haiman
showed that if $\mathcal{H}_{\mu}$ and $\mathrm{gr} R_{\mu}$ have the
same Hilbert series, then $\mathcal{H}_{\mu} = \mathrm{gr}
R_{\mu}$. While this would demonstrate the $n!$ Conjecture, the
obvious problem is that one needs first to know the Hilbert series of
$\mathcal{H}_{\mu}$, in which case the dimension can be directly
calculated. The way around this problem lies in the theory of Orbit
Harmonics developed by Garsia and Haiman. The main result is the
following.

\begin{theorem}[\cite{GHOH}]
  Let $\Phi_{\mu}$ be a basis for $R_{\mu}$.  Let $F_{\mu}(q,t) =
  \sum_{\vp \in \Phi_{\mu}} \widehat{\vp}(t,\ldots,t;q,\ldots,q)$,
  where $\widehat{\vp}$ is the leading term of $\vp$. If $F_{\mu}$ is
  symmetric in the following sense,
  \begin{equation}
    \left[ t^i q^j \right] \ F_{\mu}(q,t) \ = \ \left[ t^{n(\mu) - i}
      q^{n(\mu') - j} \right] \ F_{\mu}(q,t),
    \label{eqn:sym}
  \end{equation}
  then $\widehat{\Phi}_{\mu} = \{\widehat{\vp} \ | \ \vp \in
  \Phi_{\mu}\}$ is a basis for $\mathrm{gr}R_{\mu}$. Moreover,
  $\mathrm{gr}R_{\mu} \cong \mathcal{H}_{\mu}$ as doubly-graded $\Sn$
  modules. In particular, $\dim \mathcal{H}_{\mu} = n!$.
\label{thm:OH}
\end{theorem}

\refthm{OH} suggests the following strategy for constructing a basis
for the Garsia-Haiman module $\mathcal{H}_{\mu}$. To each filling $S$
of $\mu$, define a polynomial $\vp_S \in \mathbb{Q}[X,Y]$ so that the
evaluation matrix $(\vp_S(p_T))$ of polynomials on orbit points is
nonsingular and the corresponding degree polynomial $F_{\mu}(q,t)$ is
symmetric in the sense of \refeq{sym}. The remainder of this paper is
devoted to carrying out this strategy in the cases when $\mu$ is a two
column shape (\refsec{2col}).

%%%%%%%%%%%%%%%%%%%%%%%%%%%%%%%%%%%%%%%%%%%%%%%%%%%%%%%%%%%%
\section{Two columns}
%%%%%%%%%%%%%%%%%%%%%%%%%%%%%%%%%%%%%%%%%%%%%%%%%%%%%%%%%%%%
\label{sec:2col}

Throughout this section, we restrict our attention to partitions with
at most two columns. Following the procedure laid out in
Section~\ref{sec:mods-OH}, we will construct a basis for $R_{\mu}$
such that the degree polynomial is symmetric. Following the idea of
the kicking basis for the Garsia-Procesi modules described in
\cite{GaPr1992}, we will construct the basis together with a linear
order on fillings of $\mu$ so that the evaluation matrix has nice
triangularity properties. While the Garsia-Procesi case results in an
upper triangular matrix with nonzero diagonal entries, our matrix will
only be block triangular with respect to the largest entry.

%%%%%%%%%%%%%%%%%%%%%%%%%%%%%%%%%%%%%%%%%%%%%%%%%%%%%%%%%%%%
\subsection{The kicking tree}
%%%%%%%%%%%%%%%%%%%%%%%%%%%%%%%%%%%%%%%%%%%%%%%%%%%%%%%%%%%%
\label{sec:2-kick}%

The {\em kicking tree of $\mu$} provides a nice visualization of the
recursive construction of the proposed basis. Though proving that the
resulting collection is a basis with symmetric Hilbert series is
better done from the recursive definition, the construction is better
motivated from this viewpoint.

To construct the kicking tree, entries will be added to an empty shape
one at a time in all possible ways in some specified order, ultimately
resulting in a total ordering for the fillings. We begin by recalling
the Garsia-Procesi ordering for row-increasing tableaux
\cite{GaPr1992}.

\begin{figure}[ht]
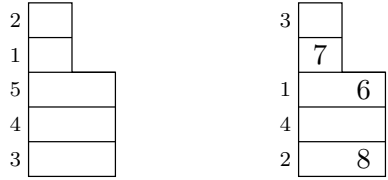

  \begin{center}
    \begin{tabular}{ccc}
      %% empty mu
      \psset{xunit=1.5em}
      \psset{yunit=1.2em}
      \psset{linewidth=.1ex}
      \pspicture(-1,0)(2,5)
      \pspolygon(0,0)(0,5)(1,5)(1,3)(2,3)(2,0)
      \psline(0,1)(2,1) \psline(0,2)(2,2) \psline(0,3)(2,3)
      \psline(0,4)(1,4) 
      \scriptsize
      \rput(-0.3,0.5){$3$}
      \rput(-0.3,1.5){$4$}
      \rput(-0.3,2.5){$5$}
      \rput(-0.3,3.5){$1$}
      \rput(-0.3,4.5){$2$}
      \normalsize 
      \endpspicture
      & \hspace{3em} &
      %% partial mu
      \psset{xunit=1.5em}
      \psset{yunit=1.2em}
      \psset{linewidth=.1ex}
      \pspicture(-1,0)(2,5)
      \pspolygon(0,0)(0,5)(1,5)(1,3)(2,3)(2,0)
      \psline(0,1)(2,1) \psline(0,2)(2,2) \psline(0,3)(2,3)
      \psline(0,4)(1,4) 
      \rput(1.5,0.5){$8$} \rput(0.5,3.5){$7$} \rput(1.5,2.5){$6$}
      \scriptsize
      \rput(-0.3,0.5){$2$}
      \rput(-0.3,1.5){$4$}
      \rput(-0.3,2.5){$1$}
      \rput(-0.3,4.5){$3$}
      \normalsize 
      \endpspicture
    \end{tabular}
    \caption{\label{fig:dominance}The row preference order for partial
      fillings.}
  \end{center}    
\end{figure}

Let $S$ be a partial filling of $\mu$ with distinct entries. Define a
total ordering on the rows of $S$ containing at least one empty cell,
called the \textit{row preference order}, as follows: empty rows of
length $2$ from top to bottom followed by (empty) rows of length $1$
from top to bottom followed by rows of length $2$ with a single
occupant beginning with the largest occupant. Given two rows $i$ and
$j$ of a (partial) filling $S$, say that \textit{k prefers row $j$
  over row $i$}, denoted $j \rowd_{k} i$, if $j$ occurs before $i$ in
the row ordering on the filling obtained by removing entries less than
$k+1$ from $S$. For example, Figure~\ref{fig:dominance} shows the
ranking of the rows (on the left) for two partial fillings of
$(2,2,2,1,1)$.

The row preference order is enough to define a total order on fillings
with unsorted rows. The basic construction of the tree is to fill
entries into unsorted rows one at a time according to row preference,
where a row of length $2$ is sorted, increasing then decreasing, as
soon as it is fully occupied. The real power of the kicking tree lies
in the weights assigned at each stage which we now describe.

Let $S$ be a partial, partially sorted filling of $\mu$ with entries
$n > n-1 > \cdots > k+1$. That is, each entry is assigned a row of
$\mu$, and an entry is assigned a specific column if and only if the
row is fully occupied. Below $S$ with arrows going down, place $k$
into a row, ordered from left to right by row preference with respect
to $k$. Label the arrow going down from $S$ to the filling with $k$ by
\begin{displaymath}
  \prod_{j \rowd_{k} \row(k)} \left( x_k - \alpha_j \right) .
\end{displaymath}
If $k$ completed a row of length $2$, say with $m>k$ already in the
row, then below this with arrows going down make two partial fillings:
the left one having $k$ before $m$ and the right having $m$ before
$k$. Label the left branch put $1$, and label the right branch
\begin{displaymath}
  \left( y_k - \beta_1 \right).
\end{displaymath}
If ignoring entries larger than $m$ does not form a rectangle, then
move the label from the arrow going down from $S$ to the left-hand
arrow just added, and add to the right-hand arrow 
\begin{displaymath}
  \prod_{\row(k) \rowd_{k} i} \hspace{-1ex} 
  \left(x_k - \alpha_i \right).
\end{displaymath}

\begin{figure}[ht]
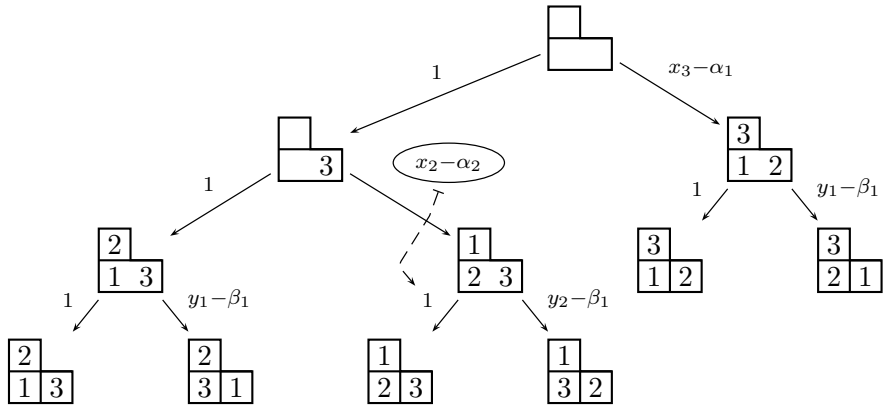

  \begin{displaymath}
    \begin{array}{ccc ccc ccc cc}
      &  &  &  &  &  &  \rnode{a1}{\tabloidXxx{}{}{}} &  &  &  \\[3ex]
      &  &  & \rnode{b1}{\tabloidXxx{}{}{3}} &  &  &  &
      & \rnode{b2}{\tabloidXxx{3}{1}{2}} & \\[3ex]
      & \rnode{c1}{\tabloidXxx{2}{1}{3}} &  &  &  
      & \rnode{c2}{\tabloidXxx{1}{2}{3}} &  &
      \rnode{c3}{\tabXxx{3}{1}{2}} & & 
      \rnode{c4}{\tabXxx{3}{2}{1}} \\[3ex]
      \rnode{d1}{\tabXxx{2}{1}{3}} & &
      \rnode{d2}{\tabXxx{2}{3}{1}} & &
      \rnode{d3}{\tabXxx{1}{2}{3}} & &
      \rnode{d4}{\tabXxx{1}{3}{2}} & & & 
    \end{array}
    \psset{nodesepA=3pt,nodesepB=3pt,linewidth=.1ex}
    \everypsbox{\scriptstyle}
    \ncline{<-}  {b1}{a1} \naput{1}
    \ncline{->}  {a1}{b2} \naput{x_3 - \alpha_1}
    \ncline{<-}  {c1}{b1} \naput{1}
    \ncline{->}  {b1}{c2} \naput{\rnode{push}{\psovalbox{x_2 - \alpha_2}}}
    \ncline{<-}  {c3}{b2} \naput{1}
    \ncline{->}  {b2}{c4} \naput{y_1 - \beta_1}
    \ncline{<-}  {d1}{c1} \naput{1}
    \ncline{->}  {c1}{d2} \naput{y_1 - \beta_1}
    \ncline{<-}  {d3}{c2} \naput{\rnode{here}{1}}
    \ncline{->}  {c2}{d4} \naput{y_2 - \beta_1}
    \ncdiag[linestyle=dashed,angleA=-110,angleB=130,linearc=.1]%
    {|->}{push}{here}
  \end{displaymath}
  \begin{center}
    \caption{\label{fig:kick21}The kicking tree for (2,1). Here the
      circled term $x_2-\alpha_2$ indicates that this term is pushed
      to the leftmost branch below. From left to right, the
      corresponding polynomials are $1 \ , \ y_1\!-\!\beta_1 \ , \
      x_2\!-\!\alpha_2 \ , \ y_2\!-\!\beta_1 \ , \ x_3\!-\!\alpha_1,
      (x_3\!-\!\alpha_1)(y_1\!-\!\beta_1)$.}
  \end{center}
\end{figure}

The tree so constructed beginning with the empty shape $\mu$ is called
the \textit{kicking tree for $\mu$}. For example, the kicking tree for
$(2,1)$ is constructed in Figure~\ref{fig:kick21}. For this example,
we omit vertical lines to indicate an unsorted row.

From the construction of the kicking tree, the product of the branch
labels from a leaf $S$ back to the empty shape $\mu$ is clearly a
polynomial. The collection of polynomials for each filling of $\mu$
forms the proposed kicking basis for $R_{\mu}$.

%%%%%%%%%%%%%%%%%%%%%%%%%%%%%%%%%%%%%%%%%%%%%%%%%%%%%%%%%%%%
\subsection{A recursive construction}
%%%%%%%%%%%%%%%%%%%%%%%%%%%%%%%%%%%%%%%%%%%%%%%%%%%%%%%%%%%%
\label{sec:2-basis}%

In order to give an alternative recursive description of the kicking
basis, we first need a bit more terminology.

\begin{figure}[ht]
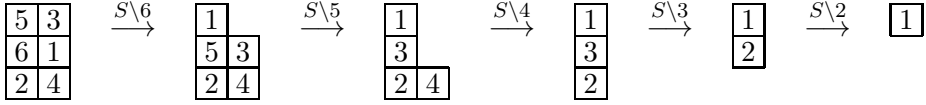

  \begin{center}
    \begin{displaymath}
      \tableau{5 & 3 \\ 6 & 1 \\ 2 & 4}\hspace{\cellsize}
      \stackrel{S\setminus 6}{\longrightarrow}\hspace{\cellsize}
      \tableau{1 \\ 5 & 3 \\ 2 & 4}\hspace{\cellsize}
      \stackrel{S\setminus 5}{\longrightarrow}\hspace{\cellsize}
      \tableau{1 \\ 3 \\ 2 & 4}\hspace{\cellsize}
      \stackrel{S\setminus 4}{\longrightarrow}\hspace{\cellsize}
      \tableau{1 \\ 3 \\ 2}\hspace{\cellsize}
      \stackrel{S\setminus 3}{\longrightarrow}\hspace{\cellsize}
      \tableau{1 \\ 2}\hspace{\cellsize}
      \stackrel{S\setminus 2}{\longrightarrow}\hspace{\cellsize}
      \tableau{1} 
    \end{displaymath}
    \caption{\label{fig:straighten}An illustration of straightening
      after removing the largest entry.}
  \end{center}
\end{figure}

For $S$ a standard filling of size $n$, define $\SN$ to be the
standard filling of size $n-1$ obtained by removing the cell
containing $n$ and {\em straightening} the shape as follows. If $n$
lies in a row of length $2$, then move the remaining cell in the same
row as $n$ above rows of length $2$ and below rows of length $1$ and
push it to the left if necessary. Otherwise slide the cells down,
preserving their order, to close the gap; see
Figure~\ref{fig:straighten}.  Notice that row dominance order commutes
with straightening.

In Definition~\ref{defn:phiS}, when the largest entry of a tableau is
removed and the remaining shape is straightened, the orbit point of
the resulting tableau is defined using the original labelling of the
rows and columns. That is, the orbit point of $\SN$ is the orbit point
of $S$ with the $n$th and $2n$th coordinates removed. For example, in
Figure~\ref{fig:straighten}, the orbit point of the filling of shape
$(2,1,1)$ will be $(\alpha_2,\alpha_1,\alpha_3,\alpha_1;
\beta_2,\beta_1,\beta_2,\beta_2)$.

\begin{definition}
  Define $\vp_{\smtableau{ _1}} = 1$. For $S$ a standard filling of $\mu$,
  $|\mu|>1$, define $\vp_S$ recursively by
  \begin{displaymath}
    \vp_S \ = \ \vp_{\SN} \cdot \! \prod_{j \rowd_{n} \row(n)} \!
    \left( x_n - \alpha_j \right) \cdot \left\{ 
      \begin{array}{cl}
        1
        & \mbox{if $n$ at the end of row$(n)$} \\[1ex]
        \left( y_k - \beta_1 \right)
        & \mbox{if $\mu=(2^b)$ and col$(n)=1$} \\[1ex]      
        \frac{\displaystyle{\left( y_k - \beta_1 \right) \hspace{-2ex}
            \prod_{\row(k) \rowd_{k} i} \hspace{-2ex} \left(x_k
              - \alpha_i \right)}}
        {\displaystyle{\prod_{j \rowd_{k} \row(k)} 
            \hspace{-1ex} \left(x_k - \alpha_j \right)}} 
        & \begin{array}{l} 
          \mbox{otherwise}
        \end{array}
      \end{array} \right.
  \end{displaymath}
  where $k$ is such that row$(k)$=row$(n)$.
\label{defn:phiS}
\end{definition}

Using the example in Figure~\ref{fig:straighten}, we compute
\begin{displaymath}
  \vp \ \raisebox{-.5\cellsize}{\smtableau{ _5 & _3 \\ _6 & _1 \\ _2 & _4}} = \
  \overbrace{(x_6 - \alpha_3) (y_1-\beta_1)}^{\tableau{6}} \cdot 
  \overbrace{(y_3-\beta_1) \left(\frac{x_3-\alpha_1}{x_3 -
        \alpha_2}\right)}^{\tableau{5}}
  \cdot \overbrace{1}^{\tableau{4}} \cdot 
  \overbrace{(x_3-\alpha_2)}^{\tableau{3}} \cdot 
  \overbrace{(x_2-\alpha_2)}^{\tableau{2}} \cdot 1,
\end{displaymath}
where each step in the recursion is indicated by the cell removed to
obtain the given terms.

The above formula associates to each standard filling $S$ of $\mu$ the
same polynomial as the kicking tree from Section~\ref{sec:2-kick}.
Notice that the denominator in the last case is precisely the label
which is `pushed down' when constructing the kicking tree. Analyzing
this statement in terms of the recursive definition yields the
following result.

\begin{proposition}
  For $S$ a standard filling of $\mu$, $\vp_S$ is a polynomial.
\end{proposition}

\begin{proof}
  The result for $\mu = (1)$ is clear, so we proceed by induction on
  $n=|\mu|$. It suffices to assume $k,n$ reside in the same length $2$
  row with $n$ in column $1$. We must show that each term in the
  denominator occurs in the numerator of $\vp_{\SN}$. The only terms
  that ever appear in any denominator are $x_i - \alpha_j$ where $i$
  lies in the second column and the entry to its left is greater. In
  particular, if $x_k - \alpha_j$ ever occurs in a numerator in the
  construction of $\vp$, it remains there through $\vp_{\SN}$. Now
  notice that the product outside of the brace (for $k$, not $n$) is
  precisely the denominator in question.
\end{proof}

To show that these polynomials form a basis for $R_{\mu}$, we show
that the evaluation matrix of polynomials on orbit points is
nonsingular. The argument uses a nested induction to show that the
matrix is almost block triangular.

\begin{theorem}
  The $n! \times n!$ matrix $\displaystyle{\left( \vp_S(p_T)
    \right)}$, where $S,T$ range over all fillings of $\mu$, is
  nonsingular. In particular, the set $\{\vp_S\}$ of polynomials
  associated to fillings of $\mu$ forms a basis for $R_{\mu}$.
  \label{thm:2-nonsing}
\end{theorem}

\begin{proof}
  We proceed by induction on $n = |\mu|$, the case $n=1$ being
  trivial. The row preference order with respect to $n$ makes
  $\displaystyle{\left( \vp_S(p_T) \right)}$ block triangular with
  respect to the row of $n$. Therefore we must show that each block,
  corresponding to $n$ in a particular row, is nonsingular. If $n$
  lies in a row of length $1$, this is immediate by induction, so
  assume $n$ lies in a row of length $2$.

  For $k < n$, let $T_k$ be a partial, partially sorted filling of
  $\mu$ with entries $k+1,k+2,\ldots,n$ (here $n$ must lie in its
  designated row of length $2$). By partially sorted, we mean that the
  row of each entry is determined, but the column is determined if and
  only if the row is fully occupied; see Figure~\ref{fig:partial} for
  an example. Let $\mathcal{T}_k$ be the set of standard fillings of
  $\mu$ which restrict to $T_k$ on $\{k+1,\ldots,n\}$, where here
  again the restriction allows the column of an entry to be
  undetermined exactly when the other occupant of the same row is at
  most $k$; again, see Figure~\ref{fig:partial}. We will show that the
  evaluation matrix for $\mathcal{T}_k$ is nonsingular by induction on
  $k$. As usual, the base case, $k=1$, is trivial.

  \begin{figure}[ht]
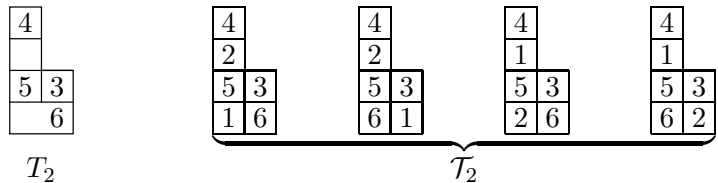

    \begin{center}
      \begin{displaymath}
        \begin{array}{ccc}
        \psset{xunit=1\cellsize}
        \psset{yunit=1\cellsize}
        \psset{linewidth=.1ex}
        \pspicture(0,0)(2,4)
        \pspolygon(0,0)(0,4)(1,4)(1,2)(2,2)(2,0)
        \psline(0,1)(2,1) \psline(0,2)(2,2) \psline(0,3)(1,3)
        \rput(0.5,3.5){$4$}
        \rput(0.5,1.5){$5$} \rput(1.5,1.5){$3$} \psline(1,1)(1,2)
        \rput(1.5,0.5){$6$}
        \endpspicture
        & \hspace{3em} &
        \underbrace{%
          \raisebox{3\cellsize}{%
          \tableau{4 \\ 2 \\5 & 3 \\ 1 & 6} \hspace{2\cellsize}
          \tableau{4 \\ 2 \\5 & 3 \\ 6 & 1} \hspace{2\cellsize}
          \tableau{4 \\ 1 \\5 & 3 \\ 2 & 6} \hspace{2\cellsize}
          \tableau{4 \\ 1 \\5 & 3 \\ 6 & 2}} } \\
        T_2 & & \mathcal{T}_2 \\[-2ex]
        \end{array}
      \end{displaymath}
      \caption{\label{fig:partial}An illustration of $T_k$ and
        $\mathcal{T}_k$.}
    \end{center}
  \end{figure}
  
  Restricting our attention to the set of polynomials and orbit points
  associated to standard fillings $S \in \mathcal{T}_k$, we put the
  following block ordering based on the position of $k$: $k$ is the
  largest entry in a row of length $2$ from highest row to lowest row;
  $k$ lies in a row of length $1$ from highest row to lowest; $k$ lies
  to the left of a larger entry from largest entry to smallest; and
  $k$ lies to the right of a larger entry again from smallest entry to
  largest. Note that the order for the first three blocks comes from
  the kicking tree, but the order of the fourth block is the reverse
  of the kicking order. By the definition of $\vp_S$, each of the four
  blocks is triangular with respect to the row of $k$, therefore by
  induction each block is nonsingular since each is a fixed
  polynomial times the polynomials associated with $\mathcal{T}_{k-1}$
  for a fixed partial filling $T_{k-1}$. 

  \begin{figure}[ht]
    \begin{displaymath}
      \left( \begin{array}{cccc}
          \raisebox{-1\cellsize}{
          \psset{xunit=.8\cellsize}
          \psset{yunit=.6\cellsize}
          \pspicture(0,0)(2,4)
          \pspolygon(0,0)(2,0)(2,2)(1,2)(1,4)(0,4)
          \psline(0,2)(2,2) 
          %\psline(1,0)(1,4)
          \rput(.9,1){$\scriptscriptstyle <  \scriptstyle k$}
          \endpspicture}
          & \ast & \rnode{l1}{\ast} & \rnode{r1}{\ast} \\
          0 &
          \raisebox{-1\cellsize}{
          \psset{xunit=.8\cellsize}
          \psset{yunit=.6\cellsize}
          \pspicture(0,0)(2,4)
          \pspolygon(0,0)(2,0)(2,2)(1,2)(1,4)(0,4)
          \psline(0,2)(2,2) 
          %\psline(1,0)(1,4)
          \rput(.5,3){$\scriptstyle k$}
          \endpspicture}
          & \rnode{l2}{\ast} & \rnode{r2}{\ast} \\
          0 & 0 &
          \rnode{l3}{%
          \raisebox{-1\cellsize}{
          \psset{xunit=.8\cellsize}
          \psset{yunit=.6\cellsize}
          \pspicture(0,0)(2,4)
          \pspolygon(0,0)(2,0)(2,2)(1,2)(1,4)(0,4)
          \psline(0,2)(2,2) 
          \psline(1,0)(1,4)
          \rput(.4,1){$\scriptstyle k$}
          \rput(1,1){$\scriptscriptstyle <$}
          \endpspicture}}
          & \rnode{r3}{\ast} \\
          \ast & \ast & 0 &
          \raisebox{-1\cellsize}{
          \psset{xunit=.8\cellsize}
          \psset{yunit=.6\cellsize}
          \pspicture(0,0)(2,4)
          \pspolygon(0,0)(2,0)(2,2)(1,2)(1,4)(0,4)
          \psline(0,2)(2,2) 
          \psline(1,0)(1,4)
          \rput(1.6,1){$\scriptstyle k$}
          \rput(1,1){$\scriptscriptstyle >$}
          \endpspicture}
        \end{array} \right)
      \psset{linewidth=.2pt, nodesep=8pt, offset=1pt}
      \ncline[nodesepA=10pt]{l1}{r1}
      \ncline[nodesepB=10pt]{r1}{l1}
      \ncline[nodesepA=10pt]{l2}{r2}
      \ncline[nodesepB=10pt]{r2}{l2}
      \ncline[nodesepA=4pt]{l3}{r3}
      \ncline[nodesepB=4pt]{r3}{l3}
    \end{displaymath}
    \caption{\label{fig:blockNR} Block structure of the evaluation
      matrix of $\mathcal{T}_{k}$.}
  \end{figure}

  Also from the definition of $\vp_S$, the first three blocks are
  triangular with respect to one another in the given order, and the
  third and fourth blocks are triangular with respect to each other as
  well. Moreover, for $S$ in one of the first three cases, the
  monomial $(y_m - \beta_i)$ does not divide $\vp_S$ for $i=1,2$ and
  any $m \geq k$ that appears by itself in a row of length two in
  $\mathcal{T}_{k}$. Therefore the block structure of the evaluation
  matrix is as depicted in Figure~\ref{fig:blockNR}. Since the first
  two blocks are nonsingular, we may perform row reductions to
  eliminate the nonzero elements in the bottom block-row of the
  matrix. These reductions will change the bottom block-row $0$ into
  some matrix, say $M$, and so by previous remarks the reductions will
  alter the fourth block by $M$ as well. Hence using the row
  reductions from the third block to restore the $0$ will also restore
  the fourth block. Hence the matrix can be made block
  triangular. Since the determinant of the matrix is the product of
  the determinants of the blocks, the full matrix is nonsingular.
\end{proof}

%%%%%%%%%%%%%%%%%%%%%%%%%%%%%%%%%%%%%%%%%%%%%%%%%%%%%%%%%%%%
\subsection{Symmetry of the Hilbert series}
%%%%%%%%%%%%%%%%%%%%%%%%%%%%%%%%%%%%%%%%%%%%%%%%%%%%%%%%%%%%
\label{sec:2-symmetry}%

Now that we have a basis for $R_{\mu}$, we must show that the
associated degree polynomial, denoted $F_{\mu}(q,t)$, is
symmetric. Recall that $F_{\mu}(q,t)$ is given by
\begin{equation}
  F_{\mu}(q,t) = \sum_{S : \mu \stackrel{\sim}{\longrightarrow} [n]}
  \widehat{\vp}_S(t,\ldots,t;q,\ldots,q) ,
\end{equation}
where $\widehat{\vp}_S$ is the highest degree term of $\vp_S$. That
is, $F_{\mu}(q,t)$ is the polynomial in $q$ and $t$ obtained by adding
leading terms of the kicking basis and recording the total $x$ degree
with $t$ and the total $y$ degree with $q$.  

Our aim is to show that $F_{\mu}(q,t)$ is symmetric, i.e.
\begin{equation}
  F_{\mu}(q,t) \ = \ t^{n(\mu)} \ q^{n(\mu')} \ F_{\mu}(1/q, 1/t).
\end{equation}
For example, from Figure~\ref{fig:kick21} we see that $F_{(2,1)}(q,t)
= 1 + 2q + 2t + qt$, which indeed exhibits the desired symmetry.

In order to establish symmetry, we will exploit a recurrence relation
that follows naturally from the recursive definition of $\vp_S$. To do
this, we must first define a more general degree polynomial, denoted
$J^{m}_{a,b}$, by
\begin{equation}
  J^{m}_{a,b}(q,t) \ = \ \frac{1}{q^m} \!\!
  \sum_{\substack{S : \mu \stackrel{\sim}{\longrightarrow} [n] \mathrm{s.t.}\\ 
      \mathrm{for} \ j = 0, \ldots, m-1 \\
      \row(n-j)=b - j, \\ \col(n-j)=1}} \!\!
  \widehat{\vp}_S (t,\ldots,t;q,\ldots,q) ,
  \label{eqn:Jabm}
\end{equation}
where $a \geq b \geq m \geq 0$. Note that $J^{m}_{a,b}$ is a
polynomial with maximum $q$ and $t$ exponents given by $b-m$ and
$\binom{a-m}{2} + \binom{b}{2}$, respectively. Pictorially,
$J^{m}_{a,b}$ is the degree polynomial of fillings of $(2^b,1^{a-b})$
with the top $m$ cells on the left-hand side of the rectangle $(2^b)$
deleted. In particular, we have
\begin{equation}
  J^{0}_{a,b}(q,t) = F_{\mu}(q,t) .
\label{eqn:JF}
\end{equation}
Therefore it is enough to show that $J^{m}_{a,b}$ is symmetric.

\begin{figure}[ht]
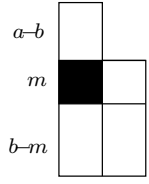

  \begin{center}
    \psset{xunit=1.5em}
    \psset{yunit=1em}
    \psset{linewidth=.1ex}
    \pspicture(0,.5)(2,6)
    \pspolygon(0,0)(0,6)(1,6)(1,4)(2,4)(2,0)
    \pspolygon[fillstyle=solid,fillcolor=black](0,4)(1,4)(1,2.5)(0,2.5)
    \psline(1,2.5)(2,2.5) \psline(1,4)(2,4) \psline(1,2.5)(1,0)
    \rput(-.7,5){$\scriptstyle a\!-\!b$}
    \rput(-.5,3.25){$\scriptstyle m$}
    \rput(-.7,1){$\scriptstyle b\!-\!m$}
    \endpspicture
    \caption{\label{fig:Jabm}The cell diagram for which $J^{m}_{a,b}$
      is the degree polynomial.}
  \end{center}
\end{figure}

\begin{proposition}
  The degree polynomials $J^{m}_{a,b}$ satisfy the following
  recurrence relations
  \begin{eqnarray}
    J^{m}_{a,b} & = & [m]_t J^{m-1}_{a-1,b-1}  +  t^{b-m}
    [a\!-\!b]_{t} J^{m}_{a-1,b}  +  t^m [b\!-\!m]_t J^{m}_{a,b-1}  + 
    q [b\!-\!m]_t J^{m+1}_{a,b}, \label{eqn:recurr1}\\
    J^{m}_{a,b} & = & t^{b-m}[m]_t J^{m-1}_{a-1,b-1} + 
    [a\!-\!b]_{t} J^{m}_{a-1,b} + q [b\!-\!m]_t J^{m}_{a,b-1} + 
    t^{a-b} [b\!-\!m]_t J^{m+1}_{a,b}, \label{eqn:recurr2}
  \end{eqnarray}
  with initial conditions
  \begin{displaymath}
    J^{b}_{a,b} = \binom{a}{b} \ [b]_{t}! \ [a-b]_{t}! 
    \hspace{2em} \mbox{and} \hspace{2em}
    J^{m}_{b,b} = J^{0}_{b,b-m},
  \end{displaymath}
  where $J^{m}_{a,b} = 0$ unless $a \geq b \geq m \geq 0$.
\label{prop:recurrence}
\end{proposition}

\begin{proof}
  The first recurrence relation follows immediately from
  Definition~\ref{defn:phiS}, with one term coming from each of the
  four regions in Figure~\ref{fig:Jabm}. The second recurrence follows
  from the first by induction on $a+b-m$, since applying the
  recurrences in either order yields the same expression.

  If all the leftmost cells in length two columns have been removed,
  then the stalk (top $a-b$ rows) and base (bottom $b$ rows) are
  independent, giving the first initial condition. The second initial
  condition is similar, since again the newly shortened rows are
  independent of the rest.
\end{proof}

The above recurrence relations follow from the recursive description
in Definition~\ref{defn:phiS}. Expanding $J^{m}_{a,b}$ twice using
both recurrence relations in Proposition~\ref{prop:recurrence} taken
in one order followed by the other establishes the desired symmetry.

\begin{theorem}
  For $a \geq b \geq m \geq 0$, we have
  \begin{displaymath}
    J^{m}_{a,b}(q,t) \ = \ t^{\binom{a-m}{2} + \binom{b}{2}} \ q^{b-m}
    \ J^{m}_{a,b}(1/q,1/t).
  \end{displaymath}
  \label{thm:symmetric}
\end{theorem}

\begin{proof}
  When $m=b$, $b=0$ or $a=b$, the result follows, so we proceed by
  induction on $a+b-m$. We will expand $J^{m}_{a,b}$ using the recurrence
  relations from Proposition~\ref{prop:recurrence}. To ease notation,
  observe that from either recurrence relation we may deduce that
  $J^{m}_{a,b}$ is divisible by $[a-b]_t! \ [b-m]_t! \
  [m]_t!$. Therefore we may define a new family polynomials, denoted
  $\wJ^{m}_{a,b}$, by
  \begin{equation}
    \wJ^{m}_{a,b} \ = \ \frac{J^{m}_{a,b}}{[a-b]_t! \ [b-m]_t! \ [m]_t!}.
  \end{equation}
  Then the recurrence relations in equations~(\ref{eqn:recurr1}) and
  (\ref{eqn:recurr2}) become
  \begin{eqnarray}
    \wJ^{m}_{a,b} & = & \wJ^{m-1}_{a-1,b-1}  +  t^{b-m}
    \wJ^{m}_{a-1,b}  +  t^m [a\!-\!b\!+\!1]_t \wJ^{m}_{a,b-1}  + 
    q [m\!+\!1]_t \wJ^{m+1}_{a,b}, \label{eqn:recurr1-mod}\\
    \wJ^{m}_{a,b} & = & t^{b-m} \wJ^{m-1}_{a-1,b-1} + 
    \wJ^{m}_{a-1,b} + q [a\!-\!b\!+\!1]_t \wJ^{m}_{a,b-1} + 
    t^{a-b} [m\!+\!1]_t \wJ^{m+1}_{a,b}, \label{eqn:recurr2-mod}
  \end{eqnarray}
  with initial conditions
  \begin{displaymath}
    J^{b}_{a,b} = \binom{a}{b}
    \hspace{2em} \mbox{and} \hspace{2em}
    J^{m}_{b,b} = J^{0}_{b,b-m}.
  \end{displaymath}
  Expanding each of the four terms in equation~(\ref{eqn:recurr1-mod})
  using the recurrence relation in equation~(\ref{eqn:recurr2-mod})
  yields the following expression for $\wJ^{m}_{a,b}$.
  \begin{displaymath}
    \begin{array}{l}
      t^{b-m} \wJ^{m-2}_{a-2,b-2} \ + \ \wJ^{m-1}_{a-2,b-1} \ + \ 
      q [a\!-\!b\!+\!1]_t \wJ^{m-1}_{a-1,b-2} \ + \ 
      t^{a-b} [m]_t \wJ^{m}_{a-1,b-1} \\[3ex]
      \ + \ t^{2(b-m)} \wJ^{m-1}_{a-2,b-1} \ + \ t^{b-m} \wJ^{m}_{a-2,b} \ + \ 
      q t^{b-m} [a\!-\!b]_t \wJ^{m}_{a-1,b-1} \ + \ 
      t^{a-m-1} [m\!+\!1]_t \wJ^{m+1}_{a-1,b} \\[3ex]
      \ + \ t^{b-1} [a\!-\!b\!+\!1]_t \wJ^{m-1}_{a-1,b-2} \ + \ 
      t^{m} [a\!-\!b\!+\!1]_t \wJ^{m}_{a-1,b-1} \ + \ 
      q t^{m} [a\!-\!b\!+\!1]_t [a\!-\!b\!+\!2]_t \wJ^{m}_{a,b-2} \\[2ex]
      \hspace{2em} \ + \ 
      t^{a-b+m-1} [a\!-\!b\!+\!1]_t [m\!+\!1]_t \wJ^{m+1}_{a,b-1} \\[3ex]
      \ + \ q t^{b-m-1} [m\!+\!1]_t \wJ^{m}_{a-1,b-1} \ + \ 
      q [m\!+\!1]_t \wJ^{m+1}_{a-1,b} \ + \ 
      q^2 [m\!+\!1]_t [a\!-\!b\!+\!1]_t \wJ^{m+1}_{a,b-1} \\[2ex]
      \hspace{2em} \ + \ 
      q t^{a-b} [m\!+\!1]_t [m\!+\!2]_t \wJ^{m+2}_{a,b}
    \end{array}
  \end{displaymath}
  Gathering terms, this expression becomes
  \begin{displaymath}
    \begin{array}{l}
      t^{b-m} \wJ^{m-2}_{a-2,b-2}
      \ + \ t^{b-m} \wJ^{m}_{a-2,b} 
      \ + \ q t^m [a\!-\!b\!+\!1]_t [a\!-\!b\!+\!2]_t \wJ^{m}_{a,b-2}
      \ + \ q t^{a-b} [m\!+\!1]_t [m\!+\!2]_t \wJ^{m+2}_{a,b} \\[3ex]
      \ + \ \left(1 + t^{2(b-m)}\right) \wJ^{m-1}_{a-2,b-1} 
      \ + \ \left(q + t^{b-1}\right) [a\!-\!b\!+\!1]_t
            \wJ^{m-1}_{a-1,b-2}
      \ + \ \left(q + t^{a-m-1}\right) [m\!+\!1]_t \wJ^{m+1}_{a-1,b}  \\[3ex]
      \ + \ \left(q^2 + t^{a-b+m-1}\right) [a\!-\!b\!+\!1]_t
            [m\!+\!1]_t \wJ^{m+1}_{a,b-1} 
      \ + \ \left(q t^{b-m-1} + t^{a-b+m}\right) \wJ^{m-1}_{a-1,b-1} \\[3ex]
      \ + \ \left(q t^{b-m} + t^{a-b}\right) [m]_t \wJ^{m-1}_{a-1,b-1}
      \ + \ \left(q t^{b-m} + t^{m}\right) [a-b]_t \wJ^{m-1}_{a-1,b-1}
    \end{array}
  \end{displaymath}
  where now each term of the above expression exhibits the desired
  symmetry, meaning the coefficient of $t^iq^j$ is equal to the
  coefficient of $t^{N-i}q^{M-j}$, where $N = (a-1)(b-m) -
  \binom{b-m}{2}$ and $M=b-m$ are the maximum powers of $t$ and $q$
  in $\wJ^{m}_{a,b}$, respectively. Hence $\wJ^{m}_{a,b}$
  symmetric. Since the product of two symmetric polynomials is
  again symmetric, the theorem now follows.
\end{proof}

In particular, by equation~(\ref{eqn:JF}), Theorem~\ref{thm:symmetric}
shows that the degree polynomial for the two column kicking basis is
indeed symmetric. Therefore by Theorem~\ref{thm:OH}, we have the
following consequence.

\begin{corollary}
  For $\mu$ a two column partition, $\{\widehat{\vp}_{S} \ | \ S:\mu
  \stackrel{\sim}{\longrightarrow} [n]\}$ is a basis for
  $\mathrm{gr}R_{\mu}$ and so too for $\mathcal{H}_{\mu}$. In
  particular, $\dim(\mathcal{H}_{\mu}) = n!$ and
  $\widetilde{K}_{\lambda,\mu}(q,t) \in \mathbb{N}[q,t]$.
\end{corollary}

%%%%%%%%%%%%%%%%%%%%%%%%%%%%%%%%%%%%%%%%%%%%%%%%%%%%%%%%%%%%
\section{Hooks}
%%%%%%%%%%%%%%%%%%%%%%%%%%%%%%%%%%%%%%%%%%%%%%%%%%%%%%%%%%%%
\label{sec:hooks}

We next treat the case of hooks, i.e. partitions $\mu =
(n-m,1^m)$. Though there exist several known bases for Garsia-Haiman
modules indexed by hooks, 
 the first in \cite{GaHa1996} and several more in
 \cite{Stembridge1994,Aval2000,Allen2002,ARR2008}.
we present this new construction because it is compatible with our two
column case, i.e. the definitions of $\vp_S$ will agree on shapes of
the form $(2,1^{n-2})$, and thus suggests how to extend this approach
to arbitrary shapes.

As with the two column case, we will construct a basis for $R_{\mu}$
such that the degree polynomial is symmetric following the idea of the
kicking basis for the Garsia-Procesi modules \cite{GaPr1992}. In this
case, the linear order on fillings of $\mu$ will have the property
that the evaluation matrix is upper triangular with nonzero diagonal
entries with respect to this basis. 

%%%%%%%%%%%%%%%%%%%%%%%%%%%%%%%%%%%%%%%%%%%%%%%%%%%%%%%%%%%%
\subsection{The kicking tree}
%%%%%%%%%%%%%%%%%%%%%%%%%%%%%%%%%%%%%%%%%%%%%%%%%%%%%%%%%%%%
\label{sec:hook-kick}%

For the case of hooks, the kicking tree is rather straightforward to
describe. Let $S$ be a partial, partially sorted filling of a hook,
say $(m,1^{M-m})$, with entries $M > M-1 > \cdots > n+1$. As before,
each entry is assigned a row, and an entry is assigned a specific
column only if its row is fully occupied. Below $S$ with arrows going
down, place $n$ into a row, ordered from left to right by row
dominance with respect to $n$. Label the arrow going down from $S$ to
the filling with $n$ by
\begin{displaymath}
  \prod_{j \rowd_{n} \row(n)} \left( x_n - \alpha_j \right) .
\end{displaymath}
If $n$ completed the bottom row, then below this will be $M$ copies of
$S \cup n$ with one entry of row 1 fixed in a column, followed by
$m-1$ copies of that with a second entry fixed in a column, and so on
for a total of $m-1$ additional levels of the tree. There are two
cases to consider in describing these levels.

First, if $n=1$ and the bottom row contains the entries $1,2, \ldots,
m$, then at each step the new level fixes the largest unfixed entry
from right to left. When $k$ is fixed into column $j$, the
preceding edge is labelled
\begin{displaymath}
  \prod_{i \in G_j} (y_k - \beta_i),
\end{displaymath}
where the product is over the set $G_j$ of unassigned column indices
greater than $j$.

Otherwise, at each step the new level fixes the smallest unfixed
letter from left to right. When $k$ is fixed into column $j$, the
preceding edge is labelled
\begin{displaymath}
  \prod_{i \in L_j} (y_k - \beta_i),
\end{displaymath}
where the product is over the set $L_j$ of unassigned column indices
less than $j$. For the smallest entry of the first row, multiply this
by the label of the arrow coming from $S$, and then replace the label
from $S$ by $1$. Finally, once the largest letter, say $M$, of the
bottom row is assigned a column, for each $h < M$ that is greater than
every other entry of the first row, remove the label from when $h$ was
added (necessarily above $S$), replacing it with $1$, and starting
from $S \cup n$, replace each $1$ with the removed label provided this
substitution was not already made at a branch above. See, for example,
\reffig{kick21}.

\begin{remark}
  At each node of the tree, the product of the labels from the empty
  shape down to that node annihilates every orbit point of a filling
  to the left of the given node. In particular, the matrix
  $\displaystyle{\left( \vp_S(p_T) \right)}$ is upper triangular with
  nonzero diagonal entries.
\label{rmk:kickh}
\end{remark}

%%%%%%%%%%%%%%%%%%%%%%%%%%%%%%%%%%%%%%%%%%%%%%%%%%%%%%%%%%%%
\subsection{A recursive basis}
%%%%%%%%%%%%%%%%%%%%%%%%%%%%%%%%%%%%%%%%%%%%%%%%%%%%%%%%%%%%
\label{sec:hook-basis}%

Both the row dominance order and straightening procedures are the same
for hooks as in the two column case, though we can describe them more
simply for this case. For row dominance, take the bottom row until
only one empty cell remains, and then take the cells from top to
bottom. For straightening, push the cells to the left and then
down. As with the two column case, when straightening a shape it is
essential to recall the original labelling of the rows and columns in
order to determine the associated orbit point.

\begin{definition}
  Define $\vp_{\smtableau{ _1}} = 1$. For $S$ a standard filling of
  $\mu$, $|\mu|=n$, define $\vp_S$ by
  \begin{displaymath}
    \vp_S \ = \vp_{\SN} \prod_{j \rowd_{n} \row(n)} \hspace{-1ex}
    \left( x_n - \alpha_j \right) \cdot \left\{ 
      \begin{array}{cr}
        \displaystyle{\prod_{\col(n) < i \leq \mu_1} \hspace{-2ex} 
          (y_n - \beta_i) }
        & \makebox[0pt][r]{if $\mu$ is a single row;} \\[1.5\vsp]
        \displaystyle{\prod_{\col(n) < \col(k)} \hspace{-2ex}
          \left(y_k - \beta_{l_k} \right) }
        & \makebox[0pt][r]{\begin{tabular}{l}%
            if $\row(n)>1$ \\ or $\col(n) > 1$,\end{tabular}} \\[1.5\vsp] 
        \displaystyle{\frac{\displaystyle{\prod_{\col(n) < \col(k)}
              \hspace{-1ex} \left(y_k - \beta_{l_k} \right)}}%
          {\displaystyle{\prod_{j \rowd_{K} 1}
              \left(x_K - \alpha_j \right)  \hspace{-2ex}
          \prod_{\substack{1 \rowd_{j} \row(j) \\ j > k \ \forall k \in
              \row(n) \setminus n}} \hspace{-3ex} \left(x_j - \alpha_1
            \right)}} }
        & \mbox{otherwise,}
      \end{array} \right.
  \end{displaymath}
  where $l_k$ is the maximum column index of all entries in row $1$
  larger than and to the left of $k$, and $K$ is the entry in the
  second column of the bottom row.
\label{defn:phiSh}
\end{definition}

For example, we compute
\begin{displaymath}
  \vp \ \raisebox{-.5\cellsize}{%
    \smtableau{ _1 \\ _4 \\ _5 & _3 & _6 & _2}} = \
  \overbrace{(y_2 - \beta_3)}^{\tableau{6}} \cdot 
  \overbrace{\frac{(y_3-\beta_1)(y_2-\beta_2)}{(x_4-\alpha_1)}}^{\tableau{5}} 
  \cdot \overbrace{(x_4-\alpha_1)(x_4-\alpha_3)}^{\tableau{4}} \cdot 
  \overbrace{\frac{(y_2-\beta_1)}{(x_2-\alpha_2)}}^{\tableau{3}} \cdot 
  \overbrace{(x_2-\alpha_2)}^{\tableau{2}} \cdot 1,
\end{displaymath}
where each step in the recursion is indicated by the cell removed to
obtain the given terms. 

Both \refprop{h-poly} and \refthm{h-nonsing} are evident from the
kicking tree description and are straightforward from the recursive
definition.

\begin{proposition}
  For $S$ a standard filling of a hook $\mu$, $\vp_S$ is a polynomial.
\label{prop:h-poly}
\end{proposition}

\begin{proof}
  The result for $\mu = (1)$ is clear, so we proceed by induction on
  $N=|\mu|$. Based on the cases above, we may assume $N$ lies in the
  first row and not in the rightmost column. Then each term appearing
  in the denominator in \refdef{phiSh} must appear in the numerator at
  some point in the construction of $\vp_{S\setminus N}$ based on the
  leftmost product. Furthermore, this is the only instance when a
  given term may appear in the denominator.
\end{proof}

\begin{theorem}
  The $n! \times n!$ evaluation matrix $\displaystyle{\left(
      \vp_S(p_T) \right)}$, where $S,T$ range over all fillings of
  $\mu$, is upper triangular with nonzero diagonal entries. In
  particular, the set $\{\vp_S\}$ forms a basis for $R_{\mu}$.
  \label{thm:h-nonsing}
\end{theorem}

\begin{proof}
  Consider the following total order on fillings of $\mu$ given by the
  kicking tree construction. Let $l$ be the largest entry such that
  $l$ lies in different rows of $S$ and $T$, or $0$ if no such entry
  exists. If $l$ is greater than some entry in the bottom row of
  either $S$ or $T$, then order $S$ and $T$ by row dominance with
  respect to $l$. Otherwise $S$ and $T$ must contain the same set of
  entries in their bottom rows and all of these are greater than
  $l$. In this case, if $l>0$, then order $S$ and $T$
  lexicographically by the bottom row, breaking a tie with row
  dominance with respect to $l$, and if $l=0$, then order $S$ and $T$
  by reverse lexicographic order. We claim that the matrix
  $\displaystyle{\left( \vp_S(p_T) \right)}$ is upper triangular with
  respect to this ordering. Since the matrix clearly has nonzero
  diagonal entries, the proposition will follow from the claim.

  We proceed by induction on $N = |\mu|$, the case $N=1$ being
  trivial. If $\mu$ is a single row, then ordering the fillings based
  on the column of $N$ gives a block triangular matrix, so by
  induction the matrix is upper triangular. In general, ordering the
  fillings based on the row of $N$ gives a block triangular matrix,
  and by induction each block corresponding to $N$ not in the first
  row may be assumed to be upper triangular with nonzero
  diagonal. Therefore it suffices to consider fillings of an honest
  hook with $N$ in the first row. 

  The result would be immediate if not for the terms in the denominator
  of $\vp_S$. The term $x_j - \alpha_1$ will evaluate to zero on any
  filling with $j$ in the bottom row. However, this term appears in
  the right hand denominator only if $y_i - \beta_1$ occurs in the
  numerator for each $i$ in the bottom row. Any filling to the left of
  this which has $j$ in the bottom row, necessarily has at least one
  of these $i$'s not in the bottom row, and so the product of the
  $y$-terms will still annihilate these fillings. Similarly, the term
  $x_k - \alpha_j$ appears in the denominator only if $y_k - \beta_1$
  appears in the numerator, so any filling annihilated by $x_k -
  \alpha_j$ must have $k$ in the first column and hence will be
  annihilated by $y_k - \beta_1$. Therefore no triangularity is lost,
  and the claim is proved.
\end{proof}

%%%%%%%%%%%%%%%%%%%%%%%%%%%%%%%%%%%%%%%%%%%%%%%%%%%%%%%%%%%%
\subsection{Symmetry of the Hilbert series}
%%%%%%%%%%%%%%%%%%%%%%%%%%%%%%%%%%%%%%%%%%%%%%%%%%%%%%%%%%%%
\label{sec:hook-symmetry}
In this section we show that the bi-graded Hilbert series of the space
spanned by $\mathcal{B}_{\mu}$ has the same symmetry as the Hilbert
series for $R_{\mu}$. 

As before, let $\Phi_{\mu}$ denote the kicking basis and define
\begin{displaymath}
  F_{\mu}(q,t) = \sum_{S : \mu \stackrel{\sim}{\longrightarrow} [n]}
  \widehat{\vp}_S (t, \ldots, t; q, \ldots, q),
\end{displaymath}
where $\widehat{\vp}_S$ is the leading term of $\vp_S$. Note that for
a hook $\mu=(n-m,1^m)$, the largest powers of $q$ and $t$ are
$(n-m)(n-m-1)/2$ and $m(m+1)/2$, which again agree with $n(\mu')$ and
$n(\mu)$, respectively.

Similar to the two column case, we can show the desired symmetry for
$F_{\mu}(q,t)$ by defining a more general function $J_{\mu}(q,t)$. By
deriving suitable recurrence relations for $F$ and $J$ in order to
establish the following theorem.

\begin{theorem}
  For $\mu$ a hook partition, both $F_{\mu}(q,t)$ and $J_{\mu}(q,t)$
  exhibit the desired symmetry. In particular, we have a basis for
  $\mathcal{H}_{\mu}$ of size $n!$, and so
  $\widetilde{K}_{\lambda,\mu}(q,t) \in \mathbb{N}[q,t]$.
  \label{thm:hook-symmetric}
\end{theorem}

%%%%%%%%%%%%%%%%%%%%%%%%%%%%%%%%%%%%%%%%%%%%%%%%%%%%%%%%%%%%
%
%  Bibliography
%
%%%%%%%%%%%%%%%%%%%%%%%%%%%%%%%%%%%%%%%%%%%%%%%%%%%%%%%%%%%%

\bibliographystyle{amsalpha}
\bibliography{nfact2.bbl}

\end{document}